\documentclass[12pt]{amsart}
\usepackage{amssymb}
\usepackage{amsfonts}
\usepackage{amsmath,amscd}
\usepackage{array}
\usepackage[matrix,frame]{xypic}
\usepackage{tableaux}

\newtheorem{example}{Example}[section]
\newtheorem{note}[example]{Note}

\newtheorem{theorem}[example]{Theorem}

\newtheorem{corollary}[example]{Corollary}

\newtheorem{proposition}[example]{Proposition}

\newtheorem{lemma}[example]{Lemma}

\def\Proof{\noindent \it Proof -- \rm}
\def\qed{\hspace{3.5mm} \hfill \vbox{\hrule height 3pt depth 2 pt width 2mm}
\bigskip}

\def\sep{\,|\,}
\def\SC{S}
\def\STab{{\rm STab}}
\def\QSym{{\it QSym}}          
\def\Sym{{\bf Sym}}            
\def\NCSF{{\bf Sym}}           
\def\FQSym{{\bf FQSym}}        
\def\FSym{{\bf FSym}}          
\def\PBT{{\bf PBT}}            

\def\PQSym{{\bf PQSym}}        
\def\WQSym{{\bf WQSym}}        

\def\Inv{{\rm Inv\,}}

\def\P{{\bf P}}
\def\M{{\bf M}}
\def\pack{{\rm pack\,}}
\def\PW{{\rm PW}}
\def\aa{{\bf a}}
\def\bb{{\bf b}}
\def\cc{{\bf c}}
\def\MM{{\mathcal M}}
\def\MC{{\mathcal M}}
\def\Ig{{\bf I}}

\newdimen\Squaresize \Squaresize=14pt
\newdimen\Thickness \Thickness=0.5pt

\def\Square#1{\hbox{\vrule width \Thickness
   \vbox to \Squaresize{\hrule height \Thickness\vss
      \hbox to \Squaresize{\hss#1\hss}
   \vss\hrule height\Thickness}
\unskip\vrule width \Thickness}
\kern-\Thickness}

\def\Vsquare#1{\vbox{\Square{$#1$}}\kern-\Thickness}

\def\std{{\rm std}}     
\def\Std{{\rm std}}     
\def\Park{{\rm park}}   

\def\<{\langle}
\def\>{\rangle}

\def\park{{\bf a}} 

\def\F{{\bf F}}         
\def\S{{\bf S}}         
\def\E{{\bf E}}         
\def\G{{\bf G}}         
\def\SG{{\mathfrak S}}  

\def\K{{\mathbb K}}

\def\T{{\mathcal T}}

\def\PF{{\rm PF}}   

\def\shuff#1#2{\mathbin{
\hbox{\vbox{ \hbox{\vrule \hskip#2 \vrule height#1 width 0pt
}%
\hrule}%
\vbox{ \hbox{\vrule \hskip#2 \vrule height#1 width 0pt
\vrule }%
\hrule}%
}}}

\long\def\psboxit#1#2{%
\begingroup\setbox0=\hbox{#2}%
\dimen0=\ht0 \advance\dimen0 by \dp0%
    \hbox{%
    \copy0%
    }
\endgroup%
}
\def\SetTableau#1#2#3#4{%
  \gdef\Tabvrule{\vrule\vrule width-0.4pt}
  \gdef\Tabhrule{\hrule\hrule height-0.4pt}  
  \gdef\Tabstrut{\vrule height#1 depth#2 width0pt\relax}
  \gdef\Tabbox##1{\hbox to #3{\hskip0.4pt\hfill\Tabstrut$#4##1$\hfill}}
} 
\def\Case#1{\vcenter{\Tabhrule%
                   \hbox{\Tabvrule\Tabbox{#1}\Tabvrule}\Tabhrule}}

\def\GenTab#1{\vcenter{\halign{&$\Case{##}$\cr#1}}\egroup}
\def\Tableau{%
  \bgroup%
  \let\ =\omit%
  \let\\=\cr%
  \offinterlineskip\GenTab}

\def\shuf{{\mathchoice{\shuff{7pt}{3.5pt}}%
{\shuff{6pt}{3pt}}%
{\shuff{4pt}{2pt}}%
{\shuff{3pt}{1.5pt}}}}%
\def\shuffle{\,\shuf\,}

\def\TD{{\mathfrak{TD}}}
\def\TC{{\mathfrak{TC}}}
\def\copNC{{ \nabla }}

\def\NI{{\rm{NI}}}
\def\PW{{\rm{PW}}}
\def\NSPW{{\rm{NSPW}}}

\title
{The $\#$ product in combinatorial Hopf algebras}

\author[J.-C. Aval, J.-C.~Novelli, and J.-Y.~Thibon]%
{Jean-Christophe Aval, Jean-Christophe Novelli and Jean-Yves Thibon}

\address[Aval]{LaBRI, Universit\'e Bordeaux I \\
351, cours de la lib\'eration \\ 33405 Talence cedex \\
FRANCE}
\address[Novelli and Thibon]{Institut Gaspard Monge, Universit\'e Paris-Est
Marne-la-Vall\'ee \\
5, Boulevard Descartes \\Champs-sur-Marne \\77454 Marne-la-Vall\'ee cedex 2 \\
FRANCE}
\email[Jean-Christophe Aval]{aval@labri.fr}
\email[Jean-Christophe Novelli]{novelli@univ-mlv.fr}
\email[Jean-Yves Thibon]{jyt@univ-mlv.fr}

\date{\today}

\begin{document}

\begin{abstract}
We show that the $\#$-product of binary trees introduced by Aval and
Viennot~[Sem. Lothar. Combin. {\bf 63} [B63h]] is in fact defined at the level
of the free associative algebra, and can be extended to most of the classical
combinatorial Hopf algebras.
\end{abstract}

\maketitle

\section{Introduction}

There is a well-known Hopf algebra structure, due to Loday and
Ronco~\cite{LR}, on the set of planar binary trees. Using a new description of
the product of this algebra, (denoted here by $\PBT$) in terms of Catalan
alternative tableaux, Aval and Viennot~\cite{AV} introduced a new product,
denoted by $\#$, which is compatible with the original graduation shifted
by~1.

Since then, Chapoton~\cite{Ch} has given a functorial interpretation
of this operation.

Most classical combinatorial Hopf algebras, including $\PBT$, admit a
realization in terms of special families of noncommutative polynomials. We
shall see that on the realization, the $\#$ product has a simple
interpretation. It can in fact be defined at the level of words over the
auxiliary alphabet. Then, it preserves in particular the algebras based on
parking functions ($\PQSym$), packed words ($\WQSym$), permutations
($\FQSym$), planar binary trees ($\PBT$), plane trees (the free tridendriform
algebra $\TD$), segmented compositions (the free cubical algebra $\TC$), Young
tableaux ($\FSym$), and integer compositions ($\NCSF$).

All definitions not recalled here can be found, \emph{e.g.},
in~\cite{NT-cras,NT-FP06,NT1}.

The present paper is an extended version of~\cite{ANT}.

\section{A semigroup of paths}

Let $A$ be an alphabet. Words over $A$ can be regarded as encoding paths in a
complete graph with a loop on each vertex, vertices being labelled by $A$.

Composition of paths, denoted by $\#$, endows the set $\Sigma(A)=A^+\cup\{0\}$
with the structure of a semigroup:
\begin{equation}
ua\# bv =
\begin{cases}
  uav &\text{if $b=a$,}\\
    0 &\text{otherwise.}
\end{cases}
\end{equation}
For example, $baaca\# adb = baacadb$ and $ab\# cd=0$. Thus, the $\#$ product
maps $A^n\times A^m$ to $A^{m+n-1}$. It is graded w.r.t. the path length
({\it i.e.}, the number of edges in the path).

We have the following obvious compatibilities with the concatenation product:
\begin{equation}
(uv)\# w = u\cdot (v\# w)\,,
\end{equation}
\begin{equation}
(u\# v)\cdot w = u\# (vw)\,.
\end{equation}

Let $d_k$ be the linear operator on $\K\<A\>$ (over some field $\K$) defined
by
\begin{equation}
d_k(w)=
\begin{cases}
  uav &\text{if $w=uaav$ for some $a$,  with $|u|=k-1$,}\\
    0 &\text{otherwise.}
\end{cases}
\end{equation}
Then, for $u$ of length $k$, 
\begin{equation}
u\# v = d_k(uv)\,.
\end{equation}

\section{Application to combinatorial Hopf algebras}

The notion of a \emph{combinatorial Hopf algebra} is a heuristic
one, referring to rich algebraic structures arising naturally
on the linear spans of various  families of combinatorial objects.
These spaces are generally endowed with several products and coproducts,
and are in particular graded connected bialgebras, hence Hopf algebras.

The most prominent combinatorial Hopf algebras can be realized in terms
of ordinary noncommutative polynomials over an auxiliary alphabet $A$.
This means that their products, which are described by combinatorial
algorithms, can be interpreted as describing the ordinary product
of certain bases of polynomials in an underlying totally
ordered alphabet  $A=\{a_1<a_2<\dots\}$.

We shall see that all these realizations are stable under the $\#$ product.
In the case of $\PBT$ (planar binary trees), we recover the result of 
Aval and Viennot~\cite{AV}.
In this case, the $\#$ product has been interpreted by Chapoton~\cite{Ch}
in representation theoretical terms.

We shall start with the most natural
algebra, $\FQSym$, based on permutations.
It contains as subalgebras  $\PBT$ (planar binary trees or the Loday-Ronco
algebra, the free dendriform algebra on one generator),
$\FSym$ (free symmetric functions, based on 
standard Young tableaux), and $\NCSF$ (noncommutative symmetric functions).

It is itself a subalgebra of
$\WQSym$, based on packed words (or set compositions),
in which the role of $\PBT$ is played by
the free dendriform trialgebra on one generator $\TD$
(based on Schr\"oder trees), the free cubical trialgebra $\TC$
(segmented compositions).

Finally, all of these algebras can be embedded in $\PQSym$, based on parking
functions.

Note that although all our algebras are actually Hopf algebras, the Hopf
structure does not play any role in this paper.

\section{Free quasi-symmetric functions: $\FQSym$ and its subalgebras}

\subsection{Free quasi-symmetric functions}

\subsubsection{The operation $d_k$ on $\FQSym$}

Recall that the alphabet $A$ is \emph{totally ordered}. Thus,
we can associate to any word over $A$ a permutation $\sigma=\std(w)$,
the \emph{standardized word} $\Std(w)$ of $w$, obtained by iteratively
scanning $w$ from left to right, and labelling $1,2,\dots$ the occurrences of
its smallest letter, then numbering the occurrences of the next one, and so
on. For example,
\begin{equation}
\Std(365182122) = 687193245.
\end{equation}

For a permutation $\sigma$, define
\begin{equation}
\G_\sigma=\sum_{\Std(w)=\sigma}w\,.
\end{equation}

We shall need the following easy property of the standardization:
\begin{lemma}
Let $=u_1u_2\dots u_n$ be a word over $A$,
and $\sigma=\sigma_1\sigma_2\dots\sigma_n = \std(u)$. Then,
for any sequence $(1\leq i_1<\dots< i_p\leq n)$,
\begin{equation}
\std(u_{i_1}u_{i_2}\dots u_{i_p})=
\std(\sigma_{i_1}\sigma_{i_2}\dots \sigma_{i_p})\,.
\end{equation}
\end{lemma}
This implies that
$\FQSym$ is stable under the $d_k$:
\begin{equation}
d_k(\G_\sigma)=
\begin{cases}
\G_{\Std(\sigma_1\dots\sigma_{k-1}\sigma_{k+1}\dots \sigma_{n})}
  & \text{if $\sigma_{k+1}=\sigma_k+1$,}\\
0 & \text{otherwise}.
\end{cases}
\end{equation}
Indeed, if $\sigma_{k+1}\not=\sigma_k+1$, no word $w$ with $w_k=w_{k+1}$ can
have $\sigma$ as standardized word, hence the $0$ answer. Otherwise, the
previous lemma applies.

We shall make use of the dual basis of the $\G_\sigma$ when dealing with
subalgebras of $\FQSym$. In the dual basis $\F_\sigma$ defined by
$\F_\sigma:=\G_{\sigma^{-1}}$, the formula is
\begin{equation}
d_k(\F_\sigma)=
\begin{cases}
\F_{\Std(\sigma_1\dots k\widehat{k+1}\dots \sigma_{n})}
  & \text{if $\sigma$ has a factor $k\,k\!+\!1$},\\
0 & \text{otherwise},
\end{cases}
\end{equation}
where $\widehat{a}$ means that $a$ is removed.

\subsubsection{Algebraic structure}

The $\G_\sigma$ span a subalgebra of the free associative algebra, denoted by
$\FQSym$.
The product is given by
\begin{equation}
\G_\alpha\G_\beta
=\sum_{\gamma=uv,\ \std(u)=\alpha,\ \std(v)=\beta}\G_\gamma\,.
\end{equation}
The set of permutations occuring in the r.h.s. is called the convolution
of $\alpha$ and $\beta$, and denoted by $\alpha*\beta$.

Hence, 
\begin{equation}
\G_\sigma \# \G_\tau
= d_k(\G_\sigma \G_\tau)
= \sum_{\nu\in \sigma\# \tau} \G_\nu,
\end{equation}
where
\begin{equation}
\sigma\# \tau = \{\nu |\, |\nu|=k+l-1, \std(\nu_1\dots \nu_k)=\sigma;
                               \std(\nu_k\dots \nu_{k+l-1}) = \tau \}\,.
\end{equation}

Indeed, $\G_\sigma\#\G_\tau$ is the sum of all words of the form $w=uxv$, with
$\std(ux)=\sigma$ and $\std(xv)=\tau$.
For example,
\begin{equation}
\G_{132} \# \G_{231} =
\G_{14352} + \G_{15342} + \G_{24351} + \G_{25341}.
\end{equation}

Note that $d_k$ induces a bijection
\begin{equation}
\SG_{n;k}:=\{\sigma\in\SG_n|\sigma_{k+1}=\sigma_k+1\}
\longrightarrow \SG_{n-1}\,.
\end{equation}
In the sequel, the notation $d_k^{-1}$ will refer to the
inverse bijection.

\subsubsection{Multiplicative bases}

The multiplicative basis $\S^\sigma$ of $\FQSym$ is defined by \cite{NCSF7}
\begin{equation}
\S^\sigma = \sum_{\tau \leq \sigma}\G_\tau\,,
\end{equation}
where $\le$ is the left weak order. 

Similarly, the multiplicative basis $\E^\sigma$ of $\FQSym$ is defined by
\cite{NCSF7}
\begin{equation}
\E^\sigma = \sum_{\tau \geq \sigma}\G_\tau\,.
\end{equation}

For $\alpha\in\SG_k$ and $\beta\in\SG_l$, define
$\alpha\vee\beta\in\SG_{k+l-1}$ as the output of
the following algorithm:

\begin{itemize}
\item scan the letters of $\alpha$ from left to right and write
\begin{equation}
\left\{
\begin{array}{ll}
\alpha_i+\beta_1-1     & \text{if } \alpha_i\leq \alpha_k,\\
\alpha_i+\max(\beta)-1 & \text{if } \alpha_i > \alpha_k,
\end{array}
\right.
\end{equation}
\item scan the letters of $\beta$ starting from the second one and write
\begin{equation}
\left\{
\begin{array}{ll}
\beta_i               & \text{if } \beta_i < \beta_1,\\
\beta_i +\alpha_k -1  & \text{if } \beta_i \geq \beta_1.
\end{array}
\right.
\end{equation}
\end{itemize}

Similarly, define $\alpha\wedge\beta$ by:

\begin{itemize}
\item scan the letters of $\alpha$ and write 
\begin{equation}
\left\{
\begin{array}{ll}
\alpha_i                & \text{if } \alpha_i < \alpha_k,\\
\alpha_i+ \beta_1 -1    & \text{if } \alpha_i \geq \alpha_k,
\end{array}
\right.
\end{equation}
\item read the letters of $\beta$ starting from the second one and write
\begin{equation}
\left\{
\begin{array}{ll}
\beta_i + \alpha_k -1    & \text{if\ } \beta_i \leq \beta_1,\\
\beta_i + \max(\alpha)-1 & \text{if\ } \beta_i > \beta_1.
\end{array}
\right.
\end{equation}
\end{itemize}

For example,
$3412\vee 35124   = 78346125$ and
$3412\wedge 35124 = 56148237$.

\begin{theorem}
\label{prodFQS}
The permutations appearing in a $\#$-product $\G_\alpha\#\G_\beta$
is an interval of the left weak order:
\begin{equation}
\G_\alpha\#\G_\beta=\sum_{\gamma\in[\alpha\wedge\beta,\alpha\vee\beta]}
\G_\gamma\,.
\end {equation}
\end{theorem}

\Proof First, it is clear that
\begin{equation}
d_k^{-1}([\alpha\wedge\beta,\alpha\vee\beta])\subseteq \alpha*\beta
\end {equation}
since for any $\gamma\in[\alpha\wedge\beta,\alpha\vee\beta]$, we have
$\std(\gamma_1\dots\gamma_k)=\alpha$ and
$\std(\gamma_{k}\dots\gamma_{k+l-1})=\beta$, so that
$d_k^{-1}(\gamma)\in\alpha *\beta$.

Let us now show the reverse inclusion, that is,
\begin{equation}
(\alpha*\beta)\cap\SG_{n;k}
\subseteq d_k^{-1}([\alpha\wedge\beta,\alpha\vee\beta]).
\end{equation}
Let $\sigma=d_k^{-1}(\alpha\vee\beta)$,
and $\tau\in(\alpha*\beta)\cap\SG_{n;k}$.
We need to show that $\Inv(\tau)\subseteq\Inv(\sigma)$.
Clearly, any inversion $(i,j)$ of $\tau$ such that $i,j\le k$ or
$i,j>k$ is also an inversion of $\sigma$, by definition of $*$.

Assume now that $i\le k<j$ and $(i,j)\not\in\Inv(\sigma)$.
By definition of $\alpha\vee\beta$, the only such pairs
$(i,j)$ are those such that $\sigma_i\le\sigma_k$ and
$\sigma_j\ge\sigma_{k+1}$.
Again, since $\std(\tau_1\dots\tau_k)=\std(\sigma_1\dots\sigma_k)$ and
$\std(\tau_k\dots\tau_{k+l-1})=\std(\sigma_k\dots\sigma_{k+l-1})$,
we have $\tau_i\le\tau_k<\tau_{k+1}\le\tau_j$, hence $(i,j)\not\in\Inv(\tau)$.

The proof of the lower bound is analogous.
\qed

Using only either the lower bound or the upper bound, one obtains:
\begin{corollary}
The bases $\S^\sigma$ and $\E^\sigma$ are multiplicative for the $\#$-product:
\begin{equation}
\S^\alpha\#\S^\beta=\S^{\alpha\vee\beta}\,
\qquad\text{and}\qquad
\E^\alpha\#\E^\beta=\E^{\alpha\wedge\beta}\,.
\end{equation}
%
\end{corollary}

For example,
\begin{equation}
\S^{3412} \# \S^{35124} = \S^{78346125}
\qquad
\E^{3412} \# \E^{35124} = \E^{56148237}.
\end{equation}

\subsubsection{Freeness}

The above description of the $\#$ product in the $\S$ basis implies the
following result:

\begin{theorem}
For the $\#$ product, $\FQSym$ is free on either $\S^\alpha$ or  $\G_\alpha$ 
where $\alpha$ runs over \emph{non-secable permutations},
that is, permutations of size $n\geq2$ such that any prefix
$\alpha_1\dots\alpha_k$ of size $2\leq k<n$ is not, up to order, the union of
an interval with maximal value $\sigma_k$ and another interval either empty or
with maximal value $n$.
\end{theorem}

\Proof
Indeed, any permutation can be uniquely decomposed as a maximal $\vee$ product
of non-secable permutations, so that the result holds for the
$\S$ basis. Now, if $\sigma=\sigma_1\vee\dots\vee\sigma_k$ is secable,
then
\begin{equation}
\G_{\sigma_1}\#\dots\# \G_{\sigma_k}=
\G_{\sigma} + \dots
\end{equation}
where the dots stand for permutations strictly smaller than $\sigma$
by Theorem~\ref{prodFQS}. Hence, $\FQSym$ is also free
on the $\G_\alpha$ with $\alpha$ non-secable.
\qed

The generating series of the number of non-secable permutations
(by shifted degree $d'(\sigma)=n-1$ for $\sigma\in\SG_n$)
is Sequence A077607 of~\cite{Slo}
\begin{equation}
\NI(t) := 2\,t + 2\,t^2 + 8\,t^3 + 44\,t^4 + 296\,t^5 + 2312\,t^6 + \dots
\end{equation}
or equivalently
\begin{equation}
1/(1-\NI(t))
=\sum_{n\geq1} n!t^{n-1}.
\end{equation}

Although there is a canonical choice of the free generators in the $\S$ basis,
there are other possibilities in the $\G$ basis, as stated by the following
proposition.

\begin{proposition}
\label{prop:NI}
For the $\#$ product, $\FQSym$ is free on the $\G_\alpha$, where $\alpha$ runs
over \emph{non-interval permutations}, that is permutations of size
$|\alpha|\ge 2$ having no prefix $\alpha_1\dots\alpha_i$ of size
$2\le i<|\alpha|$ which is up to order an interval $[j,i+j-1]$.
\end{proposition}

The starting point is the following lemma.

\begin{lemma}
\label{lemma:NI1}
Let $\sigma$ and $\tau$ be two permutations, of respective sizes $k$ and
$\ell$.
There is exactly one permutation $\gamma$ in the set $\sigma\#\tau$ such that
the set of values $\{\gamma_1,\gamma_2,\dots,\gamma_k\}$ is an interval of
integers.
We denote this permutation $\gamma$ by $\sigma\bullet\tau$.
\end{lemma}

\Proof
Let $\gamma$ be a permutation of $\sigma\#\tau$. We first observe that the
entry $\gamma_k$ is fixed and equal to 
\begin{equation}
\gamma_k = \sigma_k + \tau_1 - 1.
\end{equation}
Next, it is easy to see that there is exactly one choice for
$(\gamma_1,\dots,\gamma_{k-1})$: the first $k-1$ entries of $\sigma$ shifted
by $\gamma_k-\sigma_k$.
When this is done, the unused entries have to be ordered with respect to
$(\tau_2,\dots,\tau_l)$.
\qed

As an example we have: 
\begin{equation}
25143 \bullet 45132 = 5847\,6\,9132.
\end{equation}

Any permutation $\alpha$ may be decomposed in a unique way as a maximal $\bullet$ product. 
It is clear that permutations involved in such a maximal $\bullet$ product are
non-interval permutations.

\begin{note}
\label{remark:NI}
{\rm 
The order of "breaking" $\alpha$ in a $\bullet$ product is irrelevant.
Moreover, no new breakpoint is generated during the decomposition process.
This comes from the following observation. Let us write 
\begin{equation}
\alpha=\alpha_1\dots\alpha_r\dots\alpha_s\dots\alpha_k.
\end{equation}
First, if we assume that $\{\alpha_1,\dots,\alpha_r\}$ is an interval, we may
write $\alpha = \sigma^{(1)}\bullet\tau^{(1)}$ with $|\sigma^{(1)}|=r$.
Then $\{\tau^{(1)}_1,\dots,\tau^{(1)}_{s-r+1}\}$ is an interval if and only if
$\{\alpha_1,\dots,\alpha_s\}$ is also an interval.

Now, if we assume that $\{\alpha_1,\dots,\alpha_s\}$ is an interval, we may
write $\alpha = \sigma^{(2)}\bullet\tau^{(2)}$ with $|\sigma^{(2)}|=s$.
Then $\{\sigma^{(2)}_1,\dots,\sigma^{(2)}_{r}\}$ is an interval if and only if
$\{\alpha_1,\dots,\alpha_r\}$ is also an interval.
}
\end{note}

Let us denote by $I(\alpha)$ the number of terms in such a maximal $\bullet$
decomposition of~$\alpha$.
We observe that 
\begin{equation}
I(\sigma\bullet\tau)=I(\sigma)+I(\tau)+1
\end{equation}
and the next lemma provides a strict inequality for any other permutation
appearing in the $\#$ product of $\sigma$ and $\tau$.

\begin{lemma}
\label{lemma:NI2}
For $\sigma$ and $\tau$ two permutations, let $\gamma$ be a permutation in the
set $\sigma\#\tau$, different from $\sigma\bullet\tau$,
then $I(\gamma) < I(\sigma)+I(\tau)+1$.
\end{lemma}

We shall first introduce some notations, state and prove a sub-lemma, then
conclude the proof of Lemma \ref{lemma:NI2}.

Let $\sigma=\sigma_1\bullet\dots\bullet\sigma_r$ with $r=I(\sigma)$ and
$\tau=\tau_1\bullet\dots\bullet\tau_s$ with $s=I(\tau)$.
The size of $\sigma_i$ (resp. $\tau_i$) will be denoted by $k_i$
(resp.  $\ell_i$).
We denote by $k := 1 + \sum_{i=1}^r (k_i-1)$ the size of $\sigma$.

\begin{lemma}
\label{lemma:NI3}
With these notations, and $\gamma\neq\sigma\bullet\tau$,
if the integer $p$ is neither in the set
$\{1 + \sum_{i=1}^j (k_i-1),\ j=1,\dots,r-1\}$ ("breakpoints in $\sigma$") nor
in the set $\{k  + \sum_{i=1}^j (l_i -1),\ j=1,\dots,s-1\}$ ("breakpoints in
$\tau$"), then $\{\gamma_1,\dots,\gamma_p\}$ is {\em not} an interval.
\end{lemma}
\Proof
For $p<k$, if $\{\gamma_1,\dots,\gamma_p\}$ is an interval, then so is
$\{\sigma_1,\dots,\sigma_p\}$, which is absurd by Note \ref{remark:NI}.
For $p>k$, if $\{\gamma_1,\dots,\gamma_p\}$ is an interval, then so is
$\{\tau_1,\dots,\tau_{p-(k-1)}\}$, which is absurd by Note \ref{remark:NI}.
Finally, for $p=k$, and because $\gamma\neq\sigma\bullet\tau$, this comes from
Lemma \ref{lemma:NI1}.
\qed

\noindent {\it Proof of Lemma \ref{lemma:NI2} -- }
The proof proceeds by induction on $r+s=I(\sigma)+I(\tau)$.
If $r+s=2$, that is $\sigma$ and $\tau$ are both non-interval, the result
comes directly from Lemma~\ref{lemma:NI3}.

If $r+s>2$, we apply Lemma~\ref{lemma:NI3} to get that either $\gamma$ is
non-interval, or it may be decomposed as a $\bullet$ product, {\em but only}
on a breakpoint relative to $\sigma$ or $\tau$.
These two cases are treated in the same way, and we assume $\gamma$ may be
written as $\gamma=\gamma^{(1)}\bullet\gamma^{(2)}$ with
$\gamma^{(2)}\in (\sigma_j\bullet\dots \sigma_r )\# \tau$ for $j>1$
and we conclude by induction.
\qed

\noindent {\it Proof of Proposition \ref{prop:NI} -- }
Let $\alpha_1,\dots,\alpha_r$ be non-interval permutations.
By applying  Lemma \ref{lemma:NI2}, we get that
\begin{equation}
\G_{\alpha_1} \# \dots \# \G_{\alpha_r} =
\G_{\alpha_1\bullet\dots\bullet\alpha_r} + \sum \G_\gamma
\end{equation}
where the $\gamma$ appearing in the sum are such that
$I(\gamma)<I(\alpha_1\bullet\dots\bullet\alpha_r) = r$.
This proves the statement by triangularity.
\qed

\begin{note}
Since the inverse of a non-interval permutation is a {\em
non-internal-interval} permutation, that is a permutation $\alpha$ of size at
least $2$ such that for any $r\ge2$ no set of consecutive values
$\{\alpha_{i+1},\alpha_{i+2},\dots,\alpha_{i+r}\}$ is the set of the first $r$
integers $\{1,\dots,r\}$,  $\FQSym$ is free on  $\F_\alpha$ where $\alpha$
runs over non-internal-interval permutations.
\end{note}

\subsection{Young tableaux: $\FSym$}

The algebra $\FSym$ of free symmetric functions \cite{NCSF6} is the subalgebra
of $\FQSym$ spanned by the coplactic classes
\begin{equation}
\S_t=\sum_{Q(w)=t}w =\sum_{P(\sigma)=t}\F_\sigma
\end{equation}
where $(P,Q)$ are the $P$-symbol and $Q$-symbol defined by the
Robinson-Schensted correspondence. This algebra is isomorphic to the algebra of tableaux defined by Poirier and Reutenauer
\cite{PR}.
We shall denote by $\STab(n)$ the standard tableaux of size $n$.
We introduce the conjugate $\bar d_k$ of the operator $d_k$ defined by $\bar d_k(\sigma)=d_k(\sigma^{-1})^{-1}$
or (if we identify an injective word with its standardisation) in an equivalent manner by:
\begin{equation}
\bar d_k(\sigma)=
\begin{cases}
  \Std(ukv) &\text{if $\sigma=uk(k+1)v$,}\\
    0 &\text{otherwise.}
\end{cases}
\end{equation}
Similarly to $d_k$, $\bar d_k$ induces a bijection
\begin{equation}
\{\sigma\in\SG_n|\sigma^{-1}_{k+1}=\sigma^{-1}_k+1\}
\longrightarrow \SG_{n-1}\,.
\end{equation}
We shall use the notation $\bar d_k^{-1}$ for the inverse bijection.

Note that
\begin{equation}
S_{\smalltableau{3&4\\1&2}} = \G_{2413} + \G_{3412},
\quad\text{so that}\quad
d_1(S_{\smalltableau{3&4\\1&2}}) = \G_{312},
\end{equation}
which does not belong to $\FSym$. 
Hence $\FSym$ is not stable under the $d_k$.
However, we have:

\begin{theorem}
$\FSym$ is stable under the $\#$-product.
\end{theorem}

\Proof
For two standard tableaux $t'\in\STab(k)$ and  $t''\in\STab(\ell)$ we have
\begin{equation}
\S_{t'}\S_{t''}
= \sum_{P(\sigma')=t',\ P(\sigma'')=t''}\F_{\sigma'}\F_{\sigma''}
= \sum_{\tau\in\sigma'\shuffle\sigma''[k]}\F_\tau\,,
\end{equation}
so that $d_k(\S_{t'}\S_{t''})$ is a sum of $\F_\alpha$ without multiplicites.
We need to show that if a given $\alpha$ occurs in this sum, all the
$\alpha'$ such that $P(\alpha')=P(\alpha)$ are also present. We may assume
that $\alpha'$ differs from $\alpha$ by application of a single plactic
relation. 
We shall prove it for the relation $acb\equiv cab$, the other one being proved
in the same way.
Denote by $\sigma_1$ and $\sigma_2$ the two permutations ($\sigma_1$ of size
$k$) such that $\bar d_k^{-1}(\alpha)\in \sigma_1\shuffle\sigma_2[k]$.
Then, if $k < a$,
$\bar d_k^{-1}(\alpha')$ belongs to $\sigma_1\shuffle\sigma'_2[k]$ where $\sigma'_2$
is obtained from $\sigma_2$ by the same plactic rewriting.
Then, if $k\geq c$,
$\bar d_k^{-1}(\alpha')$ belongs to $\sigma'_1\shuffle\sigma_2[k]$ where $\sigma'_1$
is obtained from $\sigma_1$ by the same plactic rewriting.
Finally, if $a \le k<c$,
$\bar d_k^{-1}(\alpha')$ also belongs to $\sigma_1\shuffle\sigma_2[k]$ since $a$
comes from $\sigma_1$ and $c$ from $\sigma_2[k]$.
\qed

For example,
\begin{eqnarray}
\S_{\smalltableau{2 \\ 1 & 3}}\,\#\,
\S_{\smalltableau{3\\ 2\\ 1}}
&=&
\S_{\smalltableau{5\\ 4\\ 2\\ 1&3}}+
\S_{\smalltableau{5\\ 2& 4\\ 1 & 3}}
\\
\S_{\smalltableau{3\\1&2}}\,\#\,
\S_{\smalltableau{3\\2\\1}}
&=&
\S_{\smalltableau{5\\4\\3\\1&2}}
\\
\S_{\smalltableau{3\\2\\1}}\,\#\,
\S_{\smalltableau{2\\1&3}}
&=&
\S_{\smalltableau{4\\3\\2\\1&5}}
\\
\S_{\smalltableau{3\\2\\1}}\,\#\,
\S_{\smalltableau{3\\1&2}}
&=&
\S_{\smalltableau{3\\ 2&5\\  1&4}}+
\S_{\smalltableau{5\\ 3\\2\\1&4}}
\end{eqnarray}
Note that those products do not have same number of terms, so that there is no
natural definition of what would be the $\#$ product on the usual
(commutative) symmetric functions.

For $T$ an injective tableau and $S$ a subset of its entries, let us denote by $T_{|S}$ the (sub-)tableau
consisting of the restriction of $T$ to its entries in $S$.
For $T$ and $T'$ two skew-tableaux, we denote their plactic equivalence (as
for words) by $T\equiv T'$, that is we can obtain $T'$ from $T$ by playing
\emph{Jeu de Taquin}.
The $\#$ product in $\FSym$ is given by the following simple combinatorial
rule:

\begin{proposition}\label{prop:FSym}
Let $T_1$ and $T_2$ be two standard tableaux of sizes $k$ and $\ell$. Then
\begin{equation}\label{eq:FSym}
\S_{T_1}\,\#\,\S_{T_2} = \sum \S_{T}
\end{equation}
where $T$ runs over standard tableaux of size $k+\ell-1$ such that 
\begin{itemize}
\item $T_{|\{1,\dots,k\}} = T_1$;
\item $T_{|\{k,\dots,k+l-1\}} \equiv T_2$. 
\end{itemize}
\end{proposition}

\Proof
We shall use the same notation $\sigma_{|S}$ for the restriction of a
permutation $\sigma$ to a subset $S$ of its entries, and identify a
permutation with its standardisation. 

Let us first consider a tableau $T$ in the left-hand side of
Equation~(\ref{eq:FSym}) and $\gamma$ a permutation of size $k+\ell-1$ such
that $P(\gamma)=T$.
Such permutations $\gamma$ are characterized by the existence of two
permutations $\sigma$ of size $k$ and $\tau$ of size $\ell$ satisfying
\begin{itemize}
\item $P(\sigma)=T_1$;
\item $P(\tau)=T_2$;
\item $\gamma'\in\sigma\shuffle\tau[k]$;
\item $\gamma=\bar d_k(\gamma')$.
 \end{itemize}
Thus we have $P(\gamma_{|\{1,\dots,k\}})=T_1$ and
$P(\gamma_{|\{k,\dots,k+\ell-1\}})=T_2$, which implies that $T$
is in the right-hand side of (\ref{eq:FSym}).

Conversely, let us consider a tableau $T$ in the right-hand side of
Equation~(\ref{eq:FSym}) and $\gamma$ a permutation of size $k+\ell-1$ such
that $P(\gamma)=T$.
One has:  $P(\gamma_{|\{1,\dots,k\}})=T_1$ and
$P(\gamma_{|\{k,\dots,k+\ell-1\}})=T_2$.
If we write $\gamma = ukv$, then we set
$\gamma'=u'k(k+1)v'=\bar d_k^{-1}(\gamma)\in\SG_{k+\ell}$.
Then we observe that:
\begin{itemize}
\item $\gamma'=u'k(k+1)v'$,
\item $P(\gamma'_{|\{1,\dots,k\}})=T_1$,
\item $P(\gamma'_{|\{k+1,\dots,k+\ell\}})=T_2$.
\end{itemize}
This implies that $P(\gamma)=T$ is in the left-hand side of (\ref{eq:FSym}).
\qed

With this description, it is easy to compute by hand
\begin{equation}
\S_{\smalltableau{4\\ 1&2&3}} \,\#\, \S_{\smalltableau{3\\ 1&2}} = 
\S_{\smalltableau{6\\ 4&5\\ 1&2&3}}\,+\,
\S_{\smalltableau{6\\ 4\\ 1&2&3&5}}\,+\,
\S_{\smalltableau{4&6\\ 1&2&3&5}}\ .
\end{equation}

\subsection{Planar binary trees: $\PBT$}

\subsubsection{Algebraic structure}

Recall that the natural basis of $\PBT$ can be defined by
\begin{equation}
\P_T=\sum_{\T(\sigma)=T}\G_\sigma
\end{equation}
where $\T(\sigma)$ is the shape of the decreasing tree of $\sigma$.

\begin{proposition}The image of a tree by $d_k$ is either
$0$ or a single tree:
\begin{equation}
d_k(\P_T)=
\begin{cases}
\P_{T'},  \\
0,
\end{cases}
\end{equation}
according to whether $k$ is the left child of $k+1$ in the unique standard
binary search tree of shape $T$ (equivalently if the $k$-th vertex in the
infix reading of $T$ has no right child), in which case $T'$ is obtained from
$T$ by contracting this edge, the result being $0$ otherwise.
\end{proposition}

\Proof
Recall that on the basis $\F_\sigma$, any $\P_T$ is the sum of the linear
extensions of the order whose Hasse diagram is the unique binary search tree
$T'$ of shape $T$ (regarded as a poset, with the root as maximal element):
\begin{equation}
\P_T=\sum_{\sigma\in{\mathcal L}(T')}\F_\sigma\,.
\end{equation}
If the right subtree of the vertex $k$ of $T'$ is nonempty, it must contain
$k+1$ so that $k+1$ will always be before $k$ in any linear extension of $T'$,
so in this case the result is $0$.  If the right subtree of $k$ is empty,
then $k+1$ must be in the left subtree of $K$, and is actually its left son.
Then, the linear extensions of $T'$ contain permutations having the factor
$k\, k\!+\!1$.  Erasing $k+1$ in those permutations and standardizing, we get
exactly the linear extensions of the tree $T''$ obtained from $T'$ by
contracting the edge $(k,k+1)$ and standardizing.
\qed

By the above result, any product $\P_{T'}\#\P_{T''}$ is in $\PBT$. We just
need to select those linear extensions which are not annihilated by $d_k$.
Since $d_k(\F_\sigma)$ is nonzero iff $\sigma$ has (as a word) a factor
$k\,k\!+\!1$, the image under $d_k$ of the surviving linear extensions are
precisely those of the poset obtained by identifying the rightmost node of
$T'$ with the leftmost node of $T''$. Thus, $\#$ is indeed the Aval-Viennot
product.

\subsubsection{Multiplicative bases}

The multiplicative basis of initial intervals
\cite{HNT}  (corresponding to the projective
elements of~\cite{Ch}) is a subset of the $\S$ basis of $\FQSym$:
\begin{equation}
H_T = \S^\tau
\end{equation}
where $\tau$ is the maximal element of the sylvester class $T$~\cite{HNT}.
These maximal elements are the $132$-avoiding permutations. 
Hence, they are preserved by the $\#$ operation, so that we recover Chapoton's
result: the $\#$ product of two projective elements is a projective element.
One can also apply the argument the other way round: since one easily checks
that the $\#$ product of two permutations avoiding the pattern $132$ also
avoids this pattern, it is a simple proof that $\PBT$ is stable under $\#$.

As in the case of $\FQSym$, the fact that the $\S$ basis is still
multiplicative for the $\#$ product implies that the product in the $\P$ basis
is an interval in the Tamari order.

\subsection{Noncommutative symmetric functions: $\Sym$}

\subsubsection{Algebraic structure}

Recall that $\Sym$ is freely generated by the noncommutative
complete functions
\begin{equation}
S_n(A)=\sum_{i_1\le i_2\le\dots\le i_n}a_{i_1}a_{i_2}\dots a_{i_n}
=\G_{12\dots n}
\end{equation}
Here, we have obviously $S_n\# S_m=S_{n+m-1}$. This implies,
for $l(I)=r$, $I=I'i_r$ and $J=j_1J''$,
\begin{equation}
S^I\#S^J = S^{I'\cdot (i_r+j_1-1)\cdot J''}
\end{equation}
and similarly
\begin{equation}
R_I\#R_J = R_{I'\cdot (i_r+j_1-1)\cdot J''}\,.
\end{equation}
For example,
\begin{equation}
R_{1512} \# R_{43} = R_{15153}.
\end{equation}

Clearly, as a $\#$-algebra, $\NCSF^+$ is the free graded associative algebra
$\K\<x,y\>$ over the two generators
\begin{equation}
x = S_2=R_2 \qquad\qquad y = \Lambda_2=R_{11}
\end{equation}
of degree $1$, the neutral element being $S_1=R_1=\Lambda_1$.

Now, define for any composition $I=(i_0,\dots,i_r)$, the binary word
\begin{equation}
b(I) = 0^{i_0-1} 1 (0^{i_1-1})  1 \dots (0^{i_r-1}).
\end{equation}
On the binary coding of a composition $I$, one can read an expression
of $R_I$, $S^I$, and $\Lambda^I$ in terms
of $\#$-products of the generators $x,y$: 
replace the concatenation product by the $\#$-product ,
replace $0$ by respectively $x$, $x$, or $y$,
and $1$ by respectively $y$, $x+y$, or $x+y$, so that
\begin{equation}
R_I := (x^{i_0-1})^\# \# y \# (x^{i_1-1})^\# \# y \# \dots \# (x^{i_r-1})^\#,
\end{equation}
\begin{equation}
S^I := (x^{i_0-1})^\# \# (x+y) \# (x^{i_1-1})^\# \# (x+y)
       \# \dots \# (x^{i_r-1})^\#,
\end{equation}
and 
\begin{equation}
\Lambda^I := (y^{i_0-1})^\# \# (x+y) \# (y^{i_1-1})^\# \# (x+y)
             \# \dots \# (y^{i_r-1})^\#.
\end{equation}

Note that in particular, the maps sending
$S^I$ either to $R_I$ or $\Lambda^I$ are algebra automorphisms. This property
will extend to a Hopf algebra automorphism with the natural coproduct.

\subsubsection{Coproduct}

In this case, we have a natural coproduct:
the one for which $x$ and $y$ are primitive:
\begin{equation}
\copNC S_2 = S_2\otimes S_1 + S_1\otimes S_2,
\end{equation}
\begin{equation}
\copNC \Lambda_2 = \Lambda_2\otimes \Lambda_1 + \Lambda_1\otimes \Lambda_2,
\end{equation}
and, the neutral element $S_1$ is grouplike
\begin{equation}
\copNC S_1 = S_1\otimes S_1.
\end{equation}
Then,
\begin{equation}
\copNC S_n = \sum_{i=1}^n \binom{n-1}{i-1} S_i\otimes S_{n+1-i}.
\end{equation}

The coproduct of generic $S^I$, $\Lambda^I$, and $R_I$ all are the same: 
since $x$ and $y$ are primitive, $x+y$ is also primitive, so that, \emph{e.g.},
\begin{equation}
\copNC R_I = \sum_{w,w' | w\shuffle w'=b(I)} C_{w,w'}^{b(I)} R_J \otimes R_K,
\end{equation}
where $J$ (resp. $K$) are the compositions whose binary words are $w$ (resp.
$w'$), and $C_{w,w'}^{b(I)}$ is the coefficient of $b(I)$ in $w\shuffle w'$.
Another way of presenting this coproduct is as follows: given $b(I)$, choose
for each element if it appears on the left or on the right of the coproduct
(hence giving $2^{|I|-1}$ terms) and compute the corresponding products of $x$
and $y$.

Hence, the maps sending $S^I$ either to
$R_I$ or $\Lambda^I$ are Hopf algebra automorphisms.

\subsubsection{Duality: quasi-symmetric functions under $\#$}

Since $\NCSF$ is isomorphic to the Hopf algebra $\K\<x,y\>$ on two primitive
generators $x$ and $y$, 
its dual is the shuffle algebra on two generators
whose coproduct is given by deconcatenation.

Since all three bases $S$, $R$, and $\Lambda$ behave in the same way for the
Hopf structure, the same holds for their dual bases, so that the bases
$M_I$, $F_I$, and the forgotten basis of $\QSym$ have the same product and
coproduct formulas. In the basis $F_I$, this is
\begin{equation}
  F_I \# F_J = \sum_{w\in b(I)\shuffle b(J)} F_K,
\end{equation}
where $K$ is the composition such that $b(K)=w$.

For example,
\begin{equation}
  F_{3} \# F_{12} = 3\,F_{14} + 2\,F_{23} + F_{32},
\end{equation}
since
\begin{equation}
 xx \shuffle yx = 3\, yxxx + 2\,xyxx + xxyx.
\end{equation}

Note that since the product is a shuffle on words in $x$ and $y$, all
elements in a product $F_I F_J$ have same length, which is
$l(I)+l(J)-1$.

\section{Word quasi-symmetric functions: $\WQSym$ and its subalgebras}

\subsection{Word quasi-symmetric functions}

\subsubsection{Algebraic structure}

Word quasi-symmetric functions are the invariants of the quasi-symmetrizing
action of the symmetric group (in the limit of an infinite alphabet),
see, {\it e.g.}, \cite{NT-FP06}.

The \emph{packed word} $u=\pack(w)$ associated with a word $w\in A^*$ is
obtained by the following process. If $b_1<b_2<\dots <b_r$ are the letters
occuring in $w$, $u$ is the image of $w$ by the homomorphism
$b_i\mapsto a_i$.
A word $u$ is said to be \emph{packed} if $\pack(u)=u$. 
Such words can be interpreted as set compositions, or as faces of the
permutohedreon, and are sometimes called pseudo-permutations \cite{KLNPS}.

As in the case of permutations, we have:
\begin{lemma}
Let $u=u_1u_2\dots u_n$ be a word over $A$,
and $v=v_1v_2\dots v_n = \pack(u)$. Then,
for any factor of $u$,
\begin{equation}
\pack(u_iu_{i+1}\dots u_j)= \pack(v_iv_{i+1}\dots v_j)\,.
\end{equation}
\end{lemma}

The natural basis of $\WQSym$, which lifts the quasi monomial basis of $QSym$,
is labelled by packed words. 
It is defined by
\begin{equation}
\M_u=\sum_{\pack(w)=u}w\,.
\end{equation}
Note that $\WQSym$ is stable under the operators $d_k$. We have
\begin{equation}
d_k(\M_w)=
\begin{cases}
\M_{w_1\dots w_{k-1}w_{k+1}\dots w_{n}}
  & \text{if $w_k=w_{k+1}$}, \\
0 & \text{otherwise}.\end{cases}
\end{equation}
so that $\WQSym$ is stable under $\#$.
In this basis, the product is given by
\begin{equation}
\M_u\M_v =\sum_{w=u'v';\ \pack(u')=u,\,\pack(v')=v}\M_w\,.
\end{equation}
Thus, $\WQSym$ is stable under $\#$, and
\begin{equation}
\M_u \# \M_v
=d_k(\M_u\M_v)
= \sum_{w\in u\# v} \M_w,
\end{equation}
where
\begin{equation}
u\# v = \{w |\, |w|=k+l-1, \pack(w_1\dots w_k)=u;
                               \pack(w_k\dots w_{k+l-1})=v \}.
\end{equation}

For example,
\begin{equation}
\M_{121} \# \M_{12} = \M_{1212} + \M_{1213} + \M_{1312}.
\end{equation}

\subsubsection{Multiplicative bases}

Recall that there exists an order on packed words generalizing the left weak
order : it is the pseudo-permutohedron order. This order has a definition in
terms of inversions (see~\cite{KLNPS}) similar to the definition of the left
weak order.
The \emph{generalized inversion set} of a given packed word $w$ is the union
of the set of pairs $(i,j)$ such that $i<j$ and $w_i>w_j$ with coefficient one
(full inversions), and the set of pairs $(i,j)$ such that $i<j$ and $w_i=w_j$
with coefficient one half (half-inversions).

One then says that two words $u$ and $v$ satisfy $u<v$ for the
pseudo-permutohedron order iff the coefficient of any pair $(i,j)$ in $u$ is
smaller than or equal to the same coefficient in $v$.

Note that the definition of $u\vee v$ and $u\wedge v$ (see the section about
$\FQSym$) does not require $u$ and $v$ to be permutations.
One then has
\begin{theorem}
\label{prodWQS}
The words appearing in the product $\M_u\#\M_v$ is an interval of the
pseudo-permutohedron order:
\begin{equation}
\M_u\#\M_v=\sum_{w\in[u\wedge v,u\vee v]} \M_w\,.
\end {equation}
\end{theorem}

\Proof
The proof of the statement is identical to that of the corresponding property
of $\FQSym$, simply taking into account half-inversions in the picture.
\qed

The multiplicative basis $\S^u$ of $\WQSym$ is defined in \cite{NT-FP06} by
\begin{equation}
\S^u = \sum_{v \le u}\M_v\,,
\end{equation}
where $\le$ is the pseudo-permutohedron order. 

\begin{proposition}
The $\S$-basis is multiplicative for the $\#$-product:
\begin{equation}
\S^u\# \S^v = \S^{u\vee v}.
\end{equation}
\end{proposition}

Similarly, the multiplicative basis $\E^u$ of $\WQSym$ is defined in
\cite{NT-FP06} by
\begin{equation}
\E^u = \sum_{v \geq u}\M_v\,.
\end{equation}

\begin{proposition}
The $\E$-basis is multiplicative for the $\#$-product:
\begin{equation}
\E^u\# \E^v = \E^{u\wedge v}.
\end{equation}
\end{proposition}

\subsubsection{Freeness}

As in the case of $\FQSym$, we can describe a set of free generators in the
$\S$ basis for the algebra $\WQSym$.

We shall say that a packed word $u$ of size $n$ is {\em secable} if there exists a prefix
$u_1\dots u_k$ of size $2\leq k<n$ such that:
\begin{itemize}
\item the set $\{u_1,\dots, u_k\}$ is, up to order the union of an interval
with maximal value $u_k$ and another interval either empty or with
maximal value the maximal entry of the whole word $u$;
\item $\{u_1,\dots, u_k\} \cap \{u_k,\dots, u_n\} = \{u_k\}$.
\end{itemize}
 
Conversely, a packed word of size at least $2$ which is not secable will be called {\em non-secable}.

\begin{theorem}
For the $\#$ product, $\WQSym$ is free on the $\S^u$ or $\M_u$ 
where $u$ runs over non-secable packed words.
\end{theorem}
\Proof
Any packed word can be uniquely decomposed as a maximal
$\vee$ product of non-secable packed words, whence the assertion on the $\S^u$.

As in the case of $\FQSym$, the result for the $\M_u$ comes by triangularity
thanks to Theorem \ref{prodWQS}.
\qed

If a packed word $u$ is weighted by $t^{|u|-1}$, the generating series
$\PW(t)$ of  (unrestricted) packed words corresponds to Sequence A000670
of~\cite{Slo}:
\begin{equation}
\PW(t) = 1+3\,t+13\,t^2+75\,t^3+541\,t^4+4683\,t^5+47293\,t^6 + \dots
\end{equation}
The generating series $NSPW$ of non-secable packed words is related to $PW$ by
\begin{equation}
\PW(t) = 1/(1-NSPW(t))
\end{equation}
which enables us to compute $NSPW$:
\begin{equation}
\NSPW(t) = 3\,t + 4\,t^2 + 24\,t^3 + 192\,t + 1872\,t^5 + 21168\,t^6 + \dots
\end{equation}

\subsection{The free tridendriform algebra $\TD$}
\def\TT{{\mathcal T}} 

The realization of the free dendriform trialgebra given in \cite{NT-cras}
involves the following construction.
With any word $w$ of length $n$, associate a plane tree $\TT(w)$ with $n+1$
leaves, as follows: if $m=\max(w)$ and if $w$ has exactly $k-1$ occurences of
$m$, write
\begin{equation}\label{wT}
w=v_1\,m\,v_2\dots v_{k-1}\,m\,v_k\,,
\end{equation}
where the $v_i$ may be empty. Then, $\TT(w)$ is the tree obtained by
grafting the subtrees $\TT(v_1),\TT(v_2),\dots,\TT(v_k)$ (in this order)
on a common root, with the initial condition $\TT(\epsilon)=\emptyset$
for the empty word.
For example, the tree associated with $243411$ is
\begin{equation}
\vcenter{\xymatrix@C=0.5mm@R=4mm{
*{} & *{} & *{} & *{} & {}\ar@{-}[dlll]\ar@{-}[d]\ar@{-}[drrrr] \\
*{} &  {}\ar@{-}[dl]\ar@{-}[dr] & *{}
& *{4} & {}\ar@{-}[dl]\ar@{-}[dr] & *{}
& *{4} & *{} & {}\ar@{-}[dll]\ar@{-}[d]\ar@{-}[drr] \\
 {}
& *{2} &  {} &  {}
& *{3} &  {} &  {}
& *{1} &  {}
& *{1} &  {} \\
}}
\end{equation}

We shall call \emph{sectors} the zones containing numbers and say that a
sector is to the left of another sector if its number is to the left of the
other one, so that the reading of all sectors from left to right of any
$\TT(w)$ gives back $w$.

Now define a polynomial by
\begin{equation}
\MM_T := \sum_{\TT(w)=T}\MM_w\,.
\end{equation}

Then, exactly as in the case of $\PBT$, we have
\begin{theorem}
\begin{equation}
d_k (\MM_T) =
\left\{
\begin{array}{l}
\MM_{T'}, \\
0,
\end{array}
\right.
\end{equation}
depending on whether the $k$-th and $k+1$-th sectors are grafted on the same
vertex or not. In the nonzero case, $T'$ is obtained from $T$ by gluing
the $k$-th and $k\!+\!1$-th sectors.
\end{theorem}

\Proof
If the $k$-th and $k\!+\!1$-th sectors of $T$ are not grafted on the same
vertex, they cannot have the same value, so that $d_k(\M_w)$ for all $w$
satisfying $\TT(w)=T$, implying that $d_k(\MM_T)=0$.

Otherwise, all $w$ such that $\TT(w)=T$ have equal values at places $k$ and
$k+1$, so that, applying $d_k$, we find all words satisfying $\TT(w)=T'$ where
$T'$ is obtained from $T$ by forgetting the $k$-th sector.
\qed

\begin{corollary}
The operation $\#$ is internal on $\TD$.
\end{corollary}

This result can be proved in a simpler way using the multiplicative basis $\S$
of $\WQSym$.
Indeed, it has been shown in \cite{NT-FP06} that $\TD$ is the subalgebra of
$\WQSym$ generated by the $\S^w$ where $w$ runs over words avoiding the
patterns $132$ and $121$.
Since by construction the $\#$-product of two such $\S$ also avoids both
patterns $132$ and $121$, we get immediately that $\#$ is internal in $\TD$.
                                                                                  
The same holds for the multiplicative basis $\E$ of $\WQSym$.
It has also been shown in \cite{NT-FP06} that $\TD$ is the subalgebra of
$\WQSym$ generated by the $\E^w$ where $w$ runs over words avoiding the
patterns $312$ and $211$, and one easily checks that by construction, the
$\#$-product of two such $\E$ also avoids both patterns $312$ and $211$.

\subsection{The free cubical algebra $\TC$}

Define a \emph{segmented composition} as a finite sequence of integers,
separated by vertical bars or commas, \emph{e.g.}, $(2,1\sep2\sep1,2)$.
We shall associate an ordinary composition with a segmented composition by
replacing the vertical bars by commas.

There is a natural bijection between segmented compositions of sum $n$ and
sequences of length $n-1$ over three symbols $<,=,>$: start with a segmented
composition $\Ig$. If $i$ is not a descent of the underlying composition of
$\Ig$, write $<$ ; otherwise, if $i$ corresponds to a comma, write $=$ ; if
$i$ corresponds to a bar, write $>$.

Now, with each word $w$ of length $n$, associate a segmented
composition $\SC(w)$,
defined as the sequence $s_1,\dots,s_{n-1}$ where $s_i$ is the
comparison sign between $w_i$ and $w_{i+1}$.
For example, given $w=1615116244543$, one gets the sequence (and the segmented
composition):
\begin{equation}
<><>=<><=<>>  \Longleftrightarrow (2|2|1,2|2,2|1|1).
\end{equation}

Given a segmented composition $\Ig$, define
\begin{equation}
\MC_\Ig=\sum_{\SC(T)=\Ig}\MM_T\,.
\end{equation}

It has been shown in~\cite{NT-FP06} that the $\MC_\Ig$ generate a Hopf
subalgebra of $\TD$ and that their product is given by
\begin{equation}
\MC_{\Ig'} \MC_{\Ig''} = \MC_{\Ig'.\Ig''}+\MC_{\Ig',\Ig''}+\MC_{\Ig'|\Ig''}.
\end{equation}
where $\Ig'.\Ig''$ is obtained by gluing the last part of $\Ig'$ with the
first part of $\Ig''$.

As before, it is easy to see that

\begin{theorem}
\begin{equation}
d_k (\MC_\Ig) =
\left\{
\begin{array}{l}
\MC_{\Ig'}, \\
0,
\end{array}
\right.
\end{equation}
depending on whether $k$ is not or is a descent of the underlying composition
of $\Ig$.
In the nonzero case, $\Ig'$ is obtained from $\Ig$ by decreasing the entry
that corresponds to the entry containing the $k$-th cell in the corresponding
composition, that is, if $\Ig=(i_1,\dots,i_\ell)$ where the $i$ are separated
by commas or vertical bars, decreasing $i_n$ where $n$ is the smallest integer
such that $i_1+\dots+i_n>k$.
\end{theorem}

\Proof
Trivial by definition of $\MM_{\Ig}$.
\qed

For the same reason, the following result is also true:
\begin{equation}
\MC_{\Ig'} \# \MC_{\Ig''} = \MC_{\Ig'.'\Ig''},
\end{equation}
where $\Ig'.'\Ig''$ amounts to glue together the last part of $\Ig'$ with the
first part of $\Ig''$ minus one, leaving the other parts unchanged.

The fact that $\#$ is internal can also be seen in a simple way, again using
the multiplicative bases $\S$ and $\E$ of $\WQSym$.

Indeed, it has been shown in \cite{NT-FP06} that $\TC$ is the subalgebra of
$\WQSym$ generated by the $\S^w$ (resp. the $\E^w$) where $w$ runs over words
avoiding the four patterns $132$, $213$, $121$, and $212$ (resp. $312$, $231$,
$212$, and $221$).
Since the $\#$-product of two such $\S$ ($\E$) also avoids all
the given patterns, we get that $\#$ is internal in $\TC$.

As a direct consequence (as in the case of $\NCSF$), one easily
sees that $(\TC,\#)$ is the free algebra on the three generators
of (shifted) degree 1
\begin{equation}
\MC_{2} \quad
\MC_{1,1}\quad
\MC_{1\sep1}.
\end{equation}
It has therefore a canonical coproduct, for which these generators
are primitive.

\section{Parking quasi-symmetric functions: $\PQSym$}

A \emph{parking function} on $[n]=\{1,2,\dots,n\}$ is a word
$\park=a_1a_2\dots a_n$ of length $n$ on $[n]$ whose non-decreasing
rearrangement $\park^\uparrow=a'_1a'_2\dots a'_n$ satisfies $a'_i\le i$ for
all $i$. We shall denote by $\PF$ the set of parking functions.

For a word $w$ over a totally ordered alphabet in which each element
has a successor, one can define \cite{NT1} a
notion of \emph{parkized word} $\Park(w)$, a parking function which reduces to
$\Std(w)$ when $w$ is a word without repeated letters.

For $w=w_1w_2\dots w_n$ on $\{1,2,\dots\}$, we set
\begin{equation}
\label{dw}
d(w):=\min \{i | \#\{w_j\leq i\}<i \}\,.
\end{equation}
If $d(w)=n+1$, then $w$ is a parking function and the algorithm terminates,
returning~$w$. Otherwise, let $w'$ be the word obtained by decrementing all
the elements of $w$ greater than $d(w)$. Then $\Park(w):=\Park(w')$. Since
$w'$ is smaller than $w$ in the lexicographic order, the algorithm terminates
and always returns a parking function.

\smallskip
For example, let $w=(3,5,1,1,11,8,8,2)$. Then $d(w)=6$ and the word
$w'=(3,5,1,1,10,7,7,2)$.
Then $d(w')=6$ and $w''=(3,5,1,1,9,6,6,2)$. Finally, $d(w'')=8$ and
$w'''= (3,5,1,1,8,6,6,2)$, which is a parking function.
Thus, $\Park(w)=(3,5,1,1,8,6,6,2)$.

\begin{lemma}
Let $u=u_1u_2\dots u_n$ be a word over $A$,
and $\cc=c_1c_2\dots c_n = \Park(u)$. Then,
for any factor of $u$,
\begin{equation}
\Park(u_iu_{i+1}\dots u_j)= \Park(c_ic_{i+1}\dots c_j)\,.
\end{equation}
\end{lemma}

Recall from \cite{NT1}, that with a parking function $\aa$, one associates the
polynomial
\begin{equation}
\G_\aa=\sum_{\Park(w)=\aa}w\,.
\end{equation}
These polynomials form a basis of a subalgebra%
\footnote{Strictly speaking, this subalgebra is rather $\PQSym^*$, the graded
dual of the Hopf algebra $\PQSym$, but both are actually isomorphic.}
$\PQSym$ of the free associative algebra over $A$.
In this basis, the product is given by
\begin{equation}
\G_\aa\G_\bb =\sum_{\cc=uv;\ \Park(u)=\aa,\,\Park(v)=\bb}\G_\cc\,.
\end{equation}
Thus, 
\begin{equation}
\G_\aa \# \G_\bb = \sum_{\cc\in \aa\# \bb} \G_\cc,
\end{equation}
where
\begin{equation}
\aa\# \bb = \{\cc |\, |\cc|=k+l-1, \Park(\cc_1\dots \cc_k)=\aa,
                               \Park(\cc_k\dots \cc_{k+l-1})=\bb \}.
\end{equation}
Indeed, $\G_\aa\#\G_\bb$ is the sum of all words of the form $w=uxv$, with
$\Park(ux)=\aa$ and $\Park(xv)=\bb$.

Note that $\PQSym$ is not stable under the operators $d_k$.
For example, $d_1(\G_{112})$ is not in $\PQSym$.
However, let $d'_k$ be the linear operator defined by
\begin{equation}
d'_k(\G_\cc)=\begin{cases}\G_{\cc_1\dots\cc_{k-1}\cc_{k+1}\dots\cc_{n}}&
\text{if $\cc_k=\cc_{k+1}$ and
$\cc_1\dots\cc_{k-1}\cc_{k+1}\dots \cc_{n}\in \PF$,}\\
0 & \text{otherwise}.\end{cases}
\end{equation}
\begin{proposition}
Then, if $\aa$ is of length $k$,
\begin{equation}
\G_\aa\# \G_\bb =d'_k(\G_\aa \G_\bb)\,.
\end{equation}
\end{proposition}

\Proof We already know that the l.h.s. is a sum of terms
$\G_\cc$ without multiplicites. This is also true of the r.h.s.,
since $d_k'$ induces a bijection between parking functions
of length $n-1$, and parking functions $\cc$ of length $n$
such that $c_k=c_{k+1}$.
Now, if a word $u$ occurs in $\G_\aa\# \G_\bb$, its antecedent
$u'={d'}_k^{-1}(u)$ satisfies $Park(u'_1\dots u'_k)=\aa$
and $\Park(u'_{k+1}\dots u'_{k+l})=\bb$, so that
$u$ occurs in $d'_k(\G_\aa\G_\bb)$. 
For the same reason, if $u$ occurs in $d'_k(\G_\aa\G_\bb)$, then
$u$ occurs in $\G_\aa\# \G_\bb$. \qed

For example,
\begin{equation}
\G_{121} \# \G_{1141} = \G_{121161} + \G_{121151} + \G_{121141}\,,
\end{equation}
and
\begin{equation}
\begin{split}
\G_{1411} \# \G_{2124}
& =
   \G_{2722126} + \G_{2722125} + \G_{2722124} + \G_{2622127} \\
&+ \G_{2622126} + \G_{2622125} + \G_{2622124} + \G_{2522127} \\
&+ \G_{2522126} + \G_{2522125} + \G_{2522124}.
\end{split}
\end{equation}

\section{Concluding remarks}

\subsection{Dendriform structures}

For those algebras which are stable under the operators $d_k$,
and which are dendriform or tridendriform, a similar structure
can be defined for the $\#$ product, by taking the images
of the partial products by $d_k$.

\subsection{Coproducts}

Since all our $\#$-algebras are free, one may endow them with bialgebra
structures by declaring primitive any complete set of free generators.
However, no canonical choice has been found, apart from the trivial
cases of $\TC$ and $\NCSF$.

\subsection{Geometric interpretations}

In $\WQSym$, the $\M_u$ can be interpreted as the characteristic
functions of all faces of the hyperplane arrangement of type $A$.
The product of $\WQSym$ gives then the decomposition of the characteristic
function of a cartesian product of faces.
Then $d_k$ takes the intersection of this product with the hyperplane
$\alpha_k=x_k-x_{k+1}=0$.


\end{document}